\input ./style/arxiv-general.cfg
\documentclass[MSNbibl,nameyear,seceqn,dvips]{arxstspdf}
\makeatletter
   \@ifpackageloaded{graphicx}{}{\usepackage{graphicx}}
\makeatother
\usepackage{flushend}
\usepackage{stfloats}

\volume{30}
\issue{2}
\pubyear{2015}
\firstpage{184}
\lastpage{198}
\doi{10.1214/14-STS502}
\docsubty{FLA}

\makeatletter

\newcommand{\lleft}{\left}
\newcommand{\rrvert}{\vert}
\newcommand{\rright}{\right}
\newcommand{\llvert}{\vert}

\newcommand{\logit}{\operatorname{logit}}
\newcommand{\Var}{\operatorname{Var}}
\newcommand{\Cov}{\operatorname{Cov}}
\newcommand{\AR}{\operatorname{AR}}
\newcommand{\E}{\mathrm{E}}
\def\pr{p}
\newcommand{\nverts}[1][]{ n^{#1}_{V}}

\def\neff{n_{\mathrm{eff}}}
\def\nedges{s}
\def\nmutual{m}
\def\ntrties{t}

\newtheorem{them}{Theorem}[section]

\newcommand{\cI}{\mathcal{I}}

\newcommand{\citepa}[1]{\citeauthor{#1}, \citeyear{#1}}

\newcommand{\eqref}[1]{(\ref{#1})}

\def\sfrac#1#2{#1/#2}

\makeatother

\begin{document}
\begin{frontmatter}

\title{On the Question of Effective Sample Size in Network Modeling:
An Asymptotic Inquiry}
\runtitle{Effective Sample Size in Network Modeling}

\begin{aug}
\author[A]{\fnms{Pavel N.}~\snm{Krivitsky}\ead[label=e1]{pavel@uow.edu.au}\ead[label=u1,url]{http://www.krivitsky.net/research/}}
\and
\author[B]{\fnms{Eric D.}~\snm{Kolaczyk}\corref{}\ead[label=e2]{kolaczyk@bu.edu}\ead[label=u2,url]{http://math.bu.edu/people/kolaczyk/}}\vspace*{-10pt}
\runauthor{P. Krivitsky and E. Kolaczyk}

\affiliation{University of Wollongong and Boston University}

\address[A]{Pavel N. Krivitsky is Lecturer, School of Mathematics and Applied Statistics and National Institute for Applied Statistics
Research Australia (NIASRA),
University of Wollongong, NSW 2522, Australia \printead{e1,u1}.}
\address[B]{Eric D. Kolaczyk is Professor, Department of Mathematics and
Statistics, Boston University, Boston, Massachusetts 02215, USA \printead{e2,u2}.}
\end{aug}

%
\begin{abstract}
The modeling and analysis of networks and network data has seen an explosion
of interest in recent years and represents an exciting direction for potential
growth in statistics. Despite the already substantial amount of work
done in
this area to date by researchers from various disciplines, however,
there remain
many questions of a decidedly foundational nature---natural analogues
of standard
questions already posed and addressed in more classical areas of
statistics---that have yet to even be posed, much less addressed. Here we raise and consider
one such question in connection with network modeling. Specifically, we ask,
``Given an observed network, what is the sample size?'' Using simple,
illustrative
examples from the class of exponential random graph models,
we show that the answer to this question can very much depend on basic
properties of
the networks expected under the model, as the number of vertices
$n_{V}$ in the network grows.
In particular, adopting the (asymptotic) scaling of the variance of the
maximum likelihood parameter estimates as a
notion of effective sample size ($n_{\mathrm{eff}}$), we show that when modeling
the overall propensity to have ties and the propensity to reciprocate
ties, whether
the networks are sparse or not under the model (i.e., having
a constant or an increasing number of ties per vertex, respectively)
is sufficient to yield an order of magnitude difference in $n_{\mathrm{eff}}$,
from $O(n_{V})$ to $O(n^{2}_{V})$. In addition, we report
simulation study results that suggest similar properties for models for
triadic (friend-of-a-friend) effects. We then explore some practical
implications of this result,
using both simulation and data on food-sharing from Lamalera, Indonesia.
\end{abstract}

%
\begin{keyword}
\kwd{Asymptotic normality}
\kwd{consistency}
\kwd{mutuality}
\kwd{triadic closure}
\kwd{exponential-family random graph model}
\kwd{maximum likelihood}\vspace*{-8pt}
\end{keyword}
\end{frontmatter}

\section{Introduction}\label{sec1}

Since roughly the mid-1990s, the study of networks has increased
dramatically. Researchers from across the sciences---including
biology, bioinformatics, computer science, economics, engineering,
mathematics, physics, sociology and statistics---are more and more
involved with the collection and statistical analysis of data
associated with networks. As a result, statistical methods and models
are being developed in this area at a furious pace, with contributions
coming from a wide spectrum of disciplines. See, for example, the work
of \citet{jackson2008social}, \citet{kolaczyk2009statistical} and
\citet{newman2010networks} for recent overviews from the perspective
of economics, statistics and statistical physics, respectively.

A cross-sectional network is typically represented mathematically by a
graph, say, $G=(V,E)$, where $V$ is a set of $\nverts$ vertices
(commonly written $V=\{1,\ldots,\nverts\}$) and $E$ is a set of
$\llvert  E\rrvert $ ties [represented as vertex pairs $(u,v)\in E$]. Ties
can be either directed [wherein $(u,v)$ is distinct from $(v,u)$] or undirected.
Prominent examples of networks represented in this fashion include the
World Wide Web graph (with vertices representing web pages and directed
ties representing hyperlinks pointing from one page to another),
protein--protein interaction networks in biology (with vertices
representing proteins and undirected ties representing an affinity for
two proteins to bind physically) and friendship networks (with vertices
representing people and ties representing friendship nominations in a
social survey).

A great deal of attention in the literature has been focused on the
natural problem of modeling networks---of the presence and absence of
their ties in particular. There is by now a wide
variety of network models that have been proposed, ranging from models
of largely mathematical interest to models designed to be fit
statistically to data. See, for example, the sources cited above or,
for a shorter treatment, the review paper
by \citet{airoldi2009survey}. The derivation and study of network
models is a unique endeavor, due to a number of factors. First, the
defining aspect of networks is their relational nature, and hence the
task is effectively one of modeling complex dependencies among the
vertices. Second, quite often there is no convenient space associated
with the network, and so the type of distance and geometry that can be
exploited in modeling other dependent phenomena, like time series and
spatial processes are not, generally, available when modeling networks.
Finally, network problems frequently are quite large, involving
hundreds if not thousands or hundreds of thousands of vertices and
their ties. Since a network of $\nverts$ vertices can in principle
have on the order of $O(\nverts[2])$ ties, in network modeling and
analysis---particularly statistical analysis of network data---the
sheer magnitude of the network can be a critical factor in this area.

Suppose that we observe a network, in the form of a directed graph $G=(V,E)$,
where $V$ is a set of $\nverts=|V|$ vertices and $E$ is a set of
ordered vertex pairs, indicating ties. We will focus on graphs with no
self-loops: $(u,u) \notin E$ for any $u\in V$. Alternatively, we may
think of $G$ in terms of its
$\nverts\times\nverts$ adjacency matrix $Y$,
where $Y_{ij} = 1$, if $(i,j)\in E$, and $0$, otherwise, with $Y_{ii}
\equiv0$.
What is our sample size in this setting? At the opening workshop of the
recent Program on Complex Networks, held in August of 2010 at the
Statistical and Applied
Mathematical Sciences Institute (SAMSI), in North Carolina, USA, this question
in fact evoked three different responses:
\begin{itemize}[(3)]
\item[(1)] it is the number of unique entries in $Y$, that is, $\nverts
(\nverts-1)$;
\item[(2)] it is the number of vertices, that is, $\nverts$; or
\item[(3)] it is the number of networks, that is, one.
\end{itemize}
Which answer is correct? And, why should it matter?

Despite the already vast literature on network modeling, to the best of
our knowledge this question has yet to be formally posed much less answered.
Closest to doing so are, perhaps, \citet{FrSn94e} and \citet
{SnBo99n}, who offer some discussion of this issue in the context of
jackknife and bootstrap estimation of variance in network contexts.
That this should be so is particularly curious given that the analogous
questions have been asked and answered in other areas involving
dependent data.
In particular, the notion of an \emph{effective sample size} has been
found to be useful in various contexts involving dependent data, including
survey sampling, time series analysis, spatial analysis and
even genetic case--control studies (\cite{ThZw84i}; \cite{YaRe11e}). Given a sample
of size $n$ in such contexts, an effective sample size---say, $\neff
$---typically is defined in connection with the variance of an estimator of
interest. An understanding of $\neff$, as a function of $n$, can help lend
important insight into a variety of fundamental and interrelated concerns,
including the precision with which inference can be done, the amount of
information
contributed by the data toward learning a parameter(s) and, more
practically, the
resources needed for data collection.

For example, in survey sampling, where nontrivial dependencies can arise
through the use of complex sampling designs, $\neff$ generally is
taken to be
the sample size necessary under simple random sampling with replacement to
obtain a variance equal to that resulting from the actual design used
(e.g., \citepa{lavrakas2008encyclopedia}). Alternatively, consider a simple
$\AR(1)$ time series model, where $(X_t-\mu) = \phi(X_{t-1}-\mu) +
Z_t$, for $|\phi|<1$
and $Z_t$ independent and identically distributed normal random variables,
with mean zero and variance $\sigma^2$. For a sample of size $n$, the sample
mean $\bar{X}_n$, the natural and unbiased estimator of $\mu$, has a
variance that
behaves asymptotically in $n$ like $\sigma^2/[n(1-\phi)^2]$.
Contrasting this
expression with $\sigma^2/n$, corresponding to the case of independent and
identically distributed $X_t$ (i.e., equivalent to the case where $\phi
\equiv0$),
the value $\neff= n(1-\phi)^2$ is sometimes interpreted as an
effective sample size.

In these and similar contexts, it is often possible to show that
whereas nominally the relevant (asymptotic) variance scales inversely
with the
sample size $n$, under dependency a different scaling obtains,
reflecting a
combination of (a) the nominal sample size $n$ and, importantly, (b) the
dependency structure in the data. Since networks are defined by relational
data and, hence, consist of random variables that are inherently dependent,
it seems not unreasonable to hope that we might similarly gain insight
into the above
question ``What is the sample size?'' in a network setting, with the
corresponding
$\neff$ expected to be some function of the number of vertices
$\nverts$, modified
by characteristics of the network structure itself.

Following a similar practice in these other fields, therefore, we will interpret
the scaling of the asymptotic variances of maximum likelihood estimates
in a
network model as an effective sample size. In this paper we provide
some initial
insight into the question of what is the effective sample size in
network modeling,
focusing on the impact of what is arguably the most fundamental of network
characteristics---sparsity. A now commonly acknowledged characteristic
of real-world networks is that the actual number of ties tends to scale
much more like the number of vertices [i.e., $O(\nverts)$] than the
number of potential ties [i.e., $O(\nverts[2])$]. Here we demonstrate
that two very
different regimes of asymptotics, corresponding to responses 1 and 2 above,
obtain for maximum likelihood estimates in the context of a simple case of
the popular exponential random graph models, under nonsparse and sparse
variants of the models. Response 3 suggests no meaningful asymptotics
other than via independent replication. These may arise in some unexpected
settings, such as with discrete-time Markov models for evolution of
networks over time (\citepa{HaFu10d}; \citepa{KrHa14s}, e.g.); however, we
do not
explore this direction here.

We will also show that the notion of regime of asymptotics relates to
the notion of consistency, as it applies to networks. \citet
{krivitsky2011adjusting} showed, informally, that their offset model
was consistent, in the sense that if the network's asymptotic regime
agreed with the model, the coefficients of the nonoffset terms would
converge to some asymptotic value. Although the results of \citet
{ShRi13c} suggest that consistency may be meaningless for linear ERGMs
with nontrivial dependence structure, our results, both theoretical and
simulated, suggest that offsets that control the asymptotic regime of
the network model can produce consistency-like properties.

As a technical aside, we note that the exponential random graph models
we consider here are only relatively simple versions of those commonly used
in practice. We choose to work with these models because (i)~they are amenable
to relatively standard tools in producing the theoretical results we require,
while, nevertheless, (ii) they are sufficient in allowing us
to highlight in a straightforward and illustrative manner our key
finding---that the question of effective sample size in network
settings can in fact be expected to be nontrivial and that the
answer in general is likely to be subtle, depending substantially on
basic model assumptions. That such insight may be obtained already for the
simplest models in this class not only speaks to the fundamental nature of
our results, but also appears to be fortunate, in that it would appear
that theoretical analysis of the key quantity involved in our calculations
becomes decidedly more delicate when even moderately more sophisticated models
are considered. We provide further comments in this direction at the
end of
this paper.

The rest of this paper is organized as follows. Some background and
definitions are provided in Section~\ref{sec:background}. Our main
results are presented in Section~\ref{sec:results}, first for the case
where ties arise as independent coin flips; second, for the case in
which flips corresponding to ties to and from a given pair of vertices
are dependent;
and, third, for the case of triadic (friend-of-a-friend) effects, which
we study via simulation. We then illustrate some practical implications
of our results through a simulation study in Section~\ref{sec:coverage}, exploring coverage of confidence intervals associated
with our asymptotic arguments, and through application to food-sharing
networks in Section~\ref{sec:example}, where we examine the extent to
which real-world data can be found to support nonsparse versus sparse
variants of our models. Finally, some additional discussion may be
found in Section~\ref{sec:discussion}.

\section{Background}\label{sec:background}

There are many models for networks. [See \citet
{kolaczyk2009statistical}, Chapter~6, or the review paper
by \citet{airoldi2009survey}.] The class of exponential random graph
models has a history going back
roughly 30 years and is particularly popular with practitioners in
social network analysis.
This class of models
specifies that the distribution of the adjacency matrix $Y$ follow an
exponential family form, that is,
$\pr_{\theta}(Y=y) \propto\exp ( \theta^\top g(y) )$,
for vectors $\theta$ of parameters and
$g(\cdot)$ of sufficient statistics. However, despite this seemingly
appealing feature, work in the
last 10 years has shown that exponential random graph models must be
handled with some care,
as both their theoretical properties and computational tractability can
be rather sensitive to
model specification. See \citet{robins2007recent}, for example, and
\citet{chatterjee2013estimating} for a more theoretical treatment.

Here we concern ourselves only with certain examples of the simplest
type of exponential random graph models,
wherein the dyads $(Y_{ij},Y_{ji})$ and $(Y_{k\ell}, Y_{\ell k})$ are
assumed independent,
for $(i,j)\ne(k,\ell)$, and identically distributed.
These independent dyad models arguably have the smallest amount of
dependency to still be interesting as network models. A variant of the
models introduced by \citet{holland1981exponential}, they
are in fact too simple to be appropriate for modeling in most
situations of practical interest. However, they are
ideal for our purposes, as they allow us to quickly obtain nontrivial
insight into the question of effective
sample size in network modeling using relatively standard tools and arguments.

Outside of Section~\ref{sec:triadic}, the models we consider are all
variations of the form
%
\begin{eqnarray}\label{eq:bern.w.recip}
\hspace*{24pt}\pr_{\alpha, \beta}(Y=y) & =& \prod_{i<j}
\frac{\exp \{\alpha
(y_{ij}+y_{ji}) + \beta y_{ij} y_{ji} \}} {1 + 2e^{\alpha} +
e^{2\alpha+ \beta}}\nonumber
\\[-8pt]\\[-8pt]
& =& \biggl(\frac{1}{1 + 2e^{\alpha} + e^{2\alpha+\beta}} \biggr)^{\bigl({\nverts\atop  2}\bigr)}\nonumber\\
&&{} \cdot\exp \bigl\{\alpha
\nedges(y) + \beta\nmutual(y) \bigr\},\nonumber
\end{eqnarray}
with sufficient statistics
%
\begin{equation}\label{eq:suff.stats}
\nedges(y)\equiv\sum_{i<j} (y_{ij}+y_{ji})
\quad \mbox{and}\quad  \nmutual(y)\equiv\sum_{i<j}
y_{ij}y_{ji} , \hspace*{-24pt}
\end{equation}
a so-called \emph{Bernoulli model with reciprocity}. The parameter
$\alpha$ governs the propensity of pairs of vertices $i$ and $j$ to
form a tie $(i,j)$, and the parameter $\beta$ governs the tendency
toward reciprocity, forming a tie $(j,i)$ that reciprocates $(i,j)$.
This model can be motivated from the independence and homogeneity
assumptions given above by an argument analogous to that of \citet
{FrSt86m} using the Hammersley--Clifford Theorem \citep{Be74s}, with
dependence graph being $D=\{\{(i,j),(j,i)\}: (i,j)\in V^2 \land i<j\}$,
the set of cliques of $D$ being $\{\{(i,j)\}: (i,j)\in V^2 \land i\ne
j\}\cup\{\{(i,j),(j,i)\}: (i,j)\in V^2 \land i<j\}$, and simplifying
for homogeneity.

\begin{figure*}
\begin{tabular}{@{}cc@{}}

\includegraphics{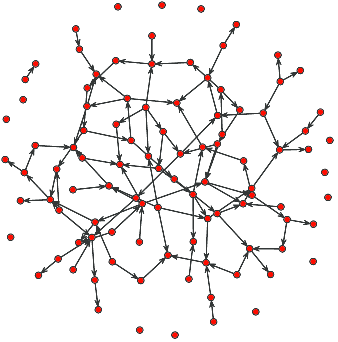}
 & \includegraphics{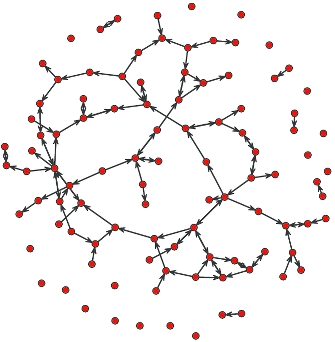}\\
\multicolumn{1}{c}{\footnotesize{(a) $\nverts=100, \nedges(y)\approx100$}}&
\multicolumn{1}{c}{\footnotesize{(b) $\nverts=100, \nedges(y)\approx100, \nmutual(y)\approx
25$}}\\[6pt]

\includegraphics{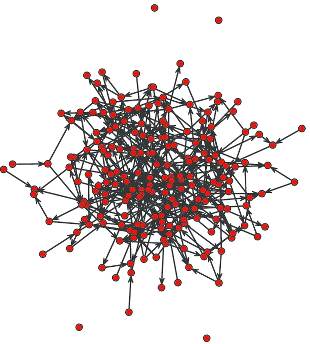}
 & \includegraphics{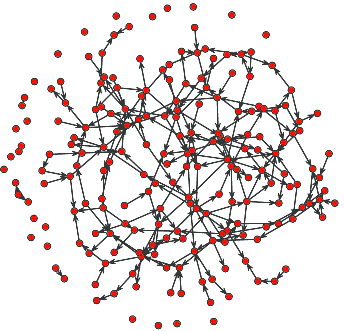}\\
\multicolumn{1}{c}{\footnotesize{(c) $\nverts=200$, preserve density of (a)}}&
\multicolumn{1}{c}{\footnotesize{(d) $\nverts=200$, preserve mean degree of (a)}}
\end{tabular}
\caption{Sampled networks drawn from four
configurations of \protect\eqref{eq:bern.w.recip}. \textup{(a)} shows a realization from a model with expected
mean degree 1 on 100 vertices, and no reciprocity effect. \textup{(b)} shows a realization from a model with
the same network size and mean degree as \textup{(a)}, but with reciprocity parameter $\beta$ set such
that the expected number of mutual ties is 25. \textup{(c)} is a realization of the model from
\textup{(a)}, scaled to 200 vertices, preserving
density; while \textup{(d)} preserves
mean degree.}\label{fig:graphs}
\end{figure*}

Of interest will be both this general model and the restricted model
$\pr_\alpha\equiv\pr_{\alpha, 0}$, wherein $\beta=0$ and there is
no reciprocity, and not just dyads, but individual potential ties
within dyads are independent. We will refer to this latter model simply
as the \emph{Bernoulli model}. Realizations of networks from this
model without and with reciprocity [holding expected tie count $\nedges
(y)$ fixed] are given in Figure~\ref{fig:graphs}(a) and (b), respectively.
%

Importantly, in both the Bernoulli model and the Bernoulli model with
reciprocity, we will examine the question of effective sample size
under both the original model parameterization and a
reparameterisation in which parameter(s) are shifted by a value $\log
\nverts$. 
\citet{krivitsky2011adjusting} introduced such shifts in an undirected
context as a way of
adjusting models like \eqref{eq:bern.w.recip} for network size such
that realizations with fixed $\alpha$ and $\beta$ would produce
network distributions with asymptotically constant expected mean
degree, $\E_{\alpha, \beta}[2\nedges(Y)/\nverts]$, for varying
$\nverts$. That
is, a configuration $(\alpha,\beta)$ that would produce a typical
$\nverts=100$ realization like that in Figure~\ref{fig:graphs}(a)
would produce an $\nverts=200$ realization like that in
Figure~\ref{fig:graphs}(d). The model's baseline asymptotic
behavior is to have a constant expected density, $\E_{\alpha,
\beta}[2\nedges(Y)/\{\nverts(\nverts-1)\}]$, such that a parameter
configuration that would produce a network
like \ref{fig:graphs}(a) for $\nverts=100$ would produce a network
like \ref{fig:graphs}(c) for $\nverts=200$.

In a directed context, ``degree'' of a given vertex $i$ is
ambiguous, as it can refer to the number of ties that vertex
makes to others ($\sum_{j\ne i}Y_{ij}$, ``outdegree''), the number of
ties others make to that vertex ($\sum_{j\ne i}Y_{ji}$,
``indegree''), the number of others to whom that vertex has at least
one connection of either type [$\sum_{j\ne i}\max(Y_{ij},Y_{ji})$],
and the number of connections that vertex has [$\sum_{j\ne
i}(Y_{ij}+Y_{ji})$]. In this work, we use either of the first
two. Then, ``mean degree'' of $Y$ is $\nedges(Y)/\nverts$, with mean
outdegree and mean indegree trivially equal; and density is
$\nedges(Y)/\{\nverts(\nverts-1)\}$.

Motivated by
similar concerns, we use the presence or absence of such shifts to
produce two different types of asymptotic behavior in our network
model classes, corresponding to sparse (asymptotically finite mean
degree) and nonsparse (asymptotically infinite mean degree) networks,
respectively. Because it is widely recognized that most large
real-world networks are sparse networks, this distinction is critical
and, as we show below, it has fundamental implications on effective
sample size and consistency.

\section{Main results}\label{sec:results}

\subsection{Bernoulli Model}
\label{sec:bernoulli}

We first present our results for the Bernoulli model. Let
$\pr_{\alpha}$ denote the model $\pr_{\alpha,0}$, as defined above, and
let $\pr^{\dag}_{\alpha}$ denote the same model, but under the mapping
$\alpha\mapsto\alpha- \log\nverts$ of the density parameter. Then,
it is easy to show that under $\pr_\alpha$ the mean vertex in- and
out-degree tends to infinity and the network density stays at
$\logit^{-1}(\alpha)$ as $\nverts\to\infty$, while under
$\pr^{\dag}_{\alpha}$, the mean degree tends to $e^{\alpha}$ while the
density tends to zero. In fact, the limiting in- and out-degree
distributions tend to a Poisson law with the stated mean.

From the perspective of traditional random graph theory, the offset model
of \citet{krivitsky2011adjusting} is asymptotically equivalent to the
standard formulation of an Erd\H{o}s--R{\'e}nyi random graph, in which the
probability of a tie scales like $e^{\alpha}/\nverts$. Alternatively,
from the
perspective of social network theory, it is useful to examine the
log-odds that
$Y_{ij}=1$, conditional on the status of all other potential ties.
Defining $Y_{[-ij]}$ to be the network $Y$ with tie $(i,j)$ removed if
present, this can be expressed as
\begin{eqnarray*}
&&\logit\pr (Y_{ij}=1 | Y_{[-ij]} = y_{[-ij]} ) \\
&&\quad \equiv\log
\frac{\pr (Y_{ij}=1 | Y_{[-ij]} = y_{[-ij]} )}{
\pr (Y_{ij}=0 | Y_{[-ij]} = y_{[-ij]} )} .
\end{eqnarray*}
This quantity goes from being a constant value $\alpha$ under
$\pr=\pr_{\alpha}$ to a value $\alpha- \log\nverts$ under $\pr
^{\dag}_{\alpha}$.
This reflects the intuition that as long as there is a cost associated
with forming and maintaining a network tie, an individual will be able
to maintain ties with a shrinking fraction of the network as the
network grows, with the average number of maintained ties being
unaffected by the growth of the network beyond a certain
point \citep{krivitsky2011adjusting}.

Given the observation of a network $Y$ randomly generated with respect
to either of these models,
initial insight into the effective sample size can be obtained by
studying the asymptotic behavior of the Fisher
information, which we denote $\cI(\alpha)$ and $\cI^{\dag}(\alpha
)$ under
$\pr_{\alpha}$ and $\pr^{\dag}_{\alpha}$, respectively.
Straightforward calculation shows that
while
\[
\cI(\alpha) = {\pmatrix{\nverts \cr 2}} \frac{2e^{\alpha
}}{(1+e^{\alpha})^2} ,
\]
in contrast,
\[
\cI^{\dag}(\alpha) = {\pmatrix{\nverts\cr 2}} \frac{2e^{\alpha
}/\nverts}{(1+e^{\alpha}/\nverts)^2} \approx
\nverts e^{\alpha} .
\]
So $\cI(\alpha) = O(\nverts[2])$, while $\cI(\alpha)^{\dag} =
O(\nverts)$, a difference by
an order of magnitude.

The implications of this difference are immediately apparent when we
consider the asymptotic behavior
of the maximum likelihood estimates of $\alpha$ under the two models.

\begin{them}\label{thm:bernoulli}
Let $\hat\alpha$ and $\hat\alpha^{\dag}$ denote the maximum
likelihood estimates of the parameter
$\alpha_0$ under models $\pr_{\alpha_0}$ and $\pr^{\dag}_{\alpha
_0}$, respectively,
where $\alpha_0\in[\alpha_{\min},\alpha_{\max}]$, for finite
$\alpha_{\min},\alpha_{\max}$.
Then under the model $\pr_{\alpha_0}$, the estimator $\hat\alpha$
is ${{\nverts\choose 2}}^{1/2}$-consistent for $\alpha_0$, and
\[
{\pmatrix{\nverts\cr 2}}^{\sfrac{1}{2}} (\hat\alpha- \alpha _0 ) \to N
\biggl( 0, \biggl\{\frac{2e^{\alpha_0}}{(1+e^{\alpha_0})^2} \biggr\} ^{-1} \biggr),
\]
while under the model $\pr^{\dag}_{\alpha_0}$, the estimator $\hat
\alpha^{\dag}$ is $\nverts[1/2]$-consistent for $\alpha_0$, and
\[
\sqrt{\nverts} \bigl(\hat\alpha^{\dag} - \alpha_0 \bigr) \to
N \bigl( 0, e^{-\alpha_0} \bigr).
\]
\end{them}

The proof of these results uses largely standard techniques for
asymptotics of estimating equations, but with a few interesting twists.
Note that, for fixed $\nverts$, the dyads $(Y_{ij},Y_{ji})$ constitute
$\nverts(\nverts-1)/2$ independent and identically distributed
bivariate random variables under both $\pr_{\alpha}$ and $\pr^{\dag
}_{\alpha}$. Consistency of the estimators in both cases can be argued
by verifying, for example, the conditions of Theorem~5.9 of \citet
{van2000asymptotic} for consistency of estimating equations. Similarly,
the proof of asymptotic normality of the estimators can be based on the
usual technique of a Taylor series expansion of the log-likelihood and,
due to the fact that we have assumed an exponential family
distribution, the asymptotic normality of the sufficient statistic
$s(y)$ in \eqref{eq:suff.stats}; however, in the case of the sparse
model $\pr_{\alpha}^{\dag}$, the dyads $\{ (Y_{ij},Y_{ji})\}_{i<j}$
follow a different distribution for each $\nverts$, and therefore an
array-based central limit theorem is required to show the asymptotic
normality of $s(y)$. But since increasing the number of vertices from,
say, $\nverts-1$ to $\nverts$, as $\nverts\rightarrow\infty$,
increases the number of dyads in our model by $\nverts-1$, a standard
triangular array central limit theorem is not appropriate here. Rather,
a double array central limit theorem is needed, such as Theorem~7.1.2
of \citet{chung2001course}. A~full derivation is provided in the
supplemental article \citep{Suppl}.

\subsection{Bernoulli Model with Reciprocity}\label
{sec:bernoulli-reciprocity}

From Theorem~\ref{thm:bernoulli} we see that
the effective sample size $\neff$
in this context can be either on the order of $\nverts$
or of $\nverts[2]$, depending on the scaling of the assumed model,
that is, on whether the model is sparse or not.
From a nonnetwork perspective, these results
can be largely anticipated by the rescaling involved, in that the
transformation $\alpha\mapsto\alpha- \log\nverts$ induces a
rescaling of the expected number of ties by $\nverts[{-1}]$.
Now, however, consider the full Bernoulli model with reciprocity, $\pr
_{\alpha, \beta}$, defined in (\ref{eq:bern.w.recip}).
Even with just two parameters the situation becomes notably more subtle.

Let $\cI(\alpha,\beta)$ be the $2\times2$ Fisher information
matrix under this model. Then calculations (not shown) completely
analogous to those
required for our previous results show that $\cI(\alpha,\beta)
= O(\nverts[2])$ and, similarly, asymptotic properties of the maximum
likelihood estimate of $(\alpha,\beta)$ analogous to those for
$\pr_\alpha$ hold.

Let us focus then on sparse versions of
$\pr_{\alpha,\beta}$. The offset used previously,
that is, mapping $\alpha$ to $\alpha- \log\nverts$, is not, by itself,
satisfactory. Call the resulting model $\pr^{\dag}_{\alpha,\beta}$.
Standard arguments show that the limiting in- and out-degree
distributions under this model will be Poisson with mean parameter
$e^{\alpha}$. On the other hand, the expected number of
\emph{reciprocated} out-ties a vertex has,
$\E^\dag_{\alpha,\beta}[2\nmutual(Y)/\nverts]$, behaves like
$e^{2\alpha+\beta}/\nverts$, and therefore tends to zero as $\nverts
\rightarrow\infty$. Thus, $\beta$ plays no role in the
limiting behavior of the model, and, indeed, reciprocity vanishes.
This fact can also be understood through examination of the Fisher
information matrix, say, $\cI^{\dag}(\alpha,\beta)$, in that direct
calculation shows that
\[
\cI^{\dag}(\alpha,\beta) = \lleft[ %
\begin{array} {c@{\quad}c} O(
\nverts) & O(1)
\\
O(1) & O(1) \end{array} %
 \rright] .
\]
That is, only the information on $\alpha$ grows with the network.
Under $p^{\dag}_{\alpha,\beta}$,
only the affinity parameter $\alpha$ can be inferred in a reliable manner.

However, the same intuition that suggests that, as the network becomes
larger, a given vertex $i$ will have an opportunity for contact with a
smaller and smaller fraction of it also suggests that if there is a
preexisting relationship in the form of a tie from $j$ to $i$, such an
opportunity likely exists for a tie from $i$ to $j$ regardless of how
large the network may
be. This, as well as direct examination of the exact expression for
the information matrix $\cI^{\dag}(\alpha,\beta)$, suggests
that the $-\log\nverts$ penalty on tie log-probability should not
apply to reciprocating ties, which may be implemented by mapping
$\beta\mapsto\beta+ \log\nverts$. Call this model, in which
$\pr^{\dag}_{\alpha,\beta}$ is augmented with this additional offset
for $\beta$, the model $\pr^{\ddag}_{\alpha,\beta}$. The
corresponding conditional log-odds of a tie now have the form
\begin{eqnarray*}
&&\logit\pr^{\ddag}_{\alpha,\beta} (Y_{ij}=1 | Y_{[-ij]} =
y_{[-ij]} )\\
&&\quad  = %
\cases{ \alpha- \log\nverts, & \mbox{if
}$y_{ji}=0$,
\cr
\alpha+ \beta, & \mbox{if }$y_{ji}=1$, }
\end{eqnarray*}
which exactly captures the intuition described.

It can be shown that under $\pr^{\ddag}_{\alpha,\beta}$ we have
$\cI^{\ddag}(\alpha,\beta) =
O(\nverts)$, indicating that information on both parameters grows at the
same rate in $\nverts$. It can also be shown that the limiting in- and
out-degree distribution is now Poisson with mean parameter $e^{\alpha}
+ e^{2\alpha+ \beta}$, and that $\E^\ddag_{\alpha,\beta
}[2\nmutual(Y)/\nverts]$ tends to
$e^{2\alpha+ \beta}$. So, both parameters play a role in the limiting
behavior of the model and the additional offset induces an
asymptotically constant expected per-vertex reciprocity in addition to
asymptotically constant expected mean degree.

Finally, we have the following analogue of Theorem~\ref{thm:bernoulli}.

\begin{them}\label{thm:mutuality}
Let $(\hat\alpha^{\ddag},\hat\beta^{\ddag})$ denote the maximum
likelihood estimate of the parameter
$(\alpha_0,\beta_0)$ under the model $\pr^{\ddag}_{\alpha_0,\beta_0}$,
where $(\alpha_0,\beta_0)\in[\alpha_{\min},\alpha_{\max}]\cdot
[\beta_{\min},\beta_{\max}]$,
for finite $\alpha_{\min},\alpha_{\max}, \beta_{\min}, \beta
_{\max}$. Then
$(\hat\alpha^{\ddag},\hat\beta^{\ddag})$ is $\nverts
[1/2]$-consistent for $(\alpha_0,\beta_0)$, and
\begin{eqnarray*}
&&\sqrt{\nverts} %
\pmatrix{ \hat\alpha^{\ddag}
- \alpha_0
\cr
\hat\beta^{\ddag} - \beta_0 } %
 \\
 &&\quad  \to N\lleft( 0 , e^{-\alpha_0} %
\lleft[\begin{array}
{c@{\quad}c} 1 & -2
\\
-2 & 4+2e^{-\alpha_0 -\beta_0} \end{array} \rright] %
 \rright) .
\end{eqnarray*}
\end{them}

Proof of this theorem, using arguments directly analogous to those of
Theorem~\ref{thm:bernoulli}, may be found in the supplemental article
\citep{Suppl}. From the theorem we see that under the sparse model
$\pr^{\ddag}_{\alpha,\beta}$, as under $\pr^{\dag}_{\alpha}$,
the effective sample size $\neff$ is $\nverts$.

\subsection{Triadic Effects}\label{sec:triadic}

Although there has been some work on obtaining closed-form
asymptotics for ERGMs with triadic---friend-of-a-friend---effects \citep{chatterjee2011random} or showing that they might not
exist \citep{ShRi13c}, these results do not appear to be directly
applicable to the per-capita asymptotic regimes that we consider in
this work. Therefore, in this section, we use simulation in an attempt to
extend the intuition developed in Section~\ref{sec:bernoulli-reciprocity}---that reciprocating ties should not be ``penalized'' for
the network size---to these triadic effects. For the sake of
simplicity, we will consider undirected networks only.

A tie between $i$ and $k$ and a tie between $k$ and $j$---that is,
that $i$ knows $k$ and $k$ knows $j$---should create a
preexisting relationship between $i$ and $j$. That is, $k$ can
``introduce'' $i$ and $j$ regardless of how large the network
is otherwise. Thus, given $i-k-j$ relationships, a potential
relationship between $i$ and $j$ should not be penalized for network
size (though $i-k-j$ themselves are); and more such two-paths
(i.e., $i-k'-j$) should have no further effect on this penalty. This
suggests an offset on the statistic called the
\emph{transitive ties} (\citepa{SnBu10i}, equation 8) or, equivalently,
Geometrically-Weighted Edgewise Shared Partners (GWESP)
\citep{morris2008specification} with parameter $\alpha$ fixed at
$0$, that is,
%
\begin{equation}
\ntrties(y)=\sum_{i<j} y_{ij} \max
_k (y_{ik},y_{jk}).\label{eq:trt}
\end{equation}
Unlike the more familiar count of the number of triangles
($\sum_{i<j<k} y_{ij} y_{ik} y_{jk}$), $\ntrties(y)$ only considers
whether a two-path between $i$ and $j$ exists, not how many of them
there are. [This
also makes it far less prone to ERGM degeneracy \citep
{schweinberger2011instability}.]

Consider the
following model, with tie count and transitive tie count
\eqref{eq:trt}:
%
\begin{eqnarray}\label{eq:bern.w.trt}
&&\pr_{\alpha, \gamma}(Y=y) \nonumber\\
&&\quad \propto \exp \bigl\{-\log(\nverts ) \bigl(
\alpha^\star\nedges(y) - \gamma^\star\ntrties(y)\bigr)\\
&&\hspace*{83pt}{} +
\alpha_0 \nedges(y) + \gamma_0 \ntrties(y) \bigr\}.\nonumber
\end{eqnarray}
As with $\pr^\dag_\alpha$, the coefficient on $\nedges$ is penalized
by network size, in the form of $\log(\nverts)\alpha^\star$, with
$\alpha^\star$ being $1$ in $\pr^\dag_\alpha$. However, the penalty
is then partially negated by increasing the coefficient on
$\ntrties$ by $\log(\nverts)\gamma^\star$. This means that, on a
sparse network,
\begin{eqnarray*}
&&\logit\pr_{\alpha,\gamma} (Y_{ij}=1 | Y_{[-ij]} =
y_{[-ij]} )
\\
&&\quad \approx\cases{ \alpha_0 - \alpha^\star\log\nverts, \quad
\mbox{if }\lnot\exists _{k\ne i,j}y_{ik}y_{kj}=1,
\cr
\alpha_0 + 3\gamma_0 - \bigl(\alpha^\star- 3
\gamma^\star\bigr)\log\nverts , \cr
\hphantom{\alpha_0 - \alpha^\star\log\nverts, \quad\, } \mbox{if }\exists_{k\ne i,j}y_{ik}y_{kj}=1.
} %
\end{eqnarray*}
This approximation holds because on an otherwise
empty network having ties $(i,k)$ and $(k,j)$, adding a tie $(i,j)$
creates not one but three transitive ties, by making all three of
the ties transitive, leading to the coefficient of $3$ on the
$\gamma$'s. However, as the network becomes more dense, this ceases to
hold exactly, because $(i,k)$
and/or $(k,j)$ may already be transitive when $(i,j)$ is added, so
only two or one transitive tie might be created.

Therefore, on a sufficiently sparse network (i.e.,
sufficiently large $\nverts$ for a given mean degree), in order to
cancel the network size penalty for a tie $(i,j)$ but retain it for
$(i,k)$ and $(k,j)$, $\gamma^\star=\alpha^\star/3$. With
$\alpha^\star=1$ per the same reasoning as before, this means that
our heuristic suggests that $\gamma^\star\approx1/3$. We verify
this empirically as follows. Define $t'(y)$---per-capita
transitive ties---as
\[
t'(y)=\frac{1}{\nverts}\sum_i\sum
_{j\ne i}y_{ij}\max_{k\ne
i,j}(y_{ik}y_{jk}).
\]
In other words, for each vertex $i$, the number of its neighbors
who have ties to at least one other neighbor of $i$ is counted, and
the resulting measures averaged over all vertices in the network. It
can be shown easily that $t'(y)\equiv2\ntrties(y)/\nverts$.

We can then ask if there exist constant values of
$\alpha^\star$ and $\gamma^\star$ that produce stable mean degree
($2\nedges(y)/\nverts$) and stable per-capita transitive ties
($t'(y)$). We constructed a series of 40 networks, sized
from 100 to 12\mbox{,}000, logarithmically spaced, for each of three
configurations of mean degree, 6, 9 and 12; and two levels of
per-capita transitivity for each: $1/2$ of the mean degree and $1/4$
of the mean degree. (This was done because per-capita transitivity
cannot exceed the mean degree.) For each combination of $\nverts$,
$2\nedges/\nverts$ and $t'/(2\nedges/\nverts)$, we used
simulated annealing to construct a network $y$ with these
statistics, and then fit an ERGM $\pr_{\alpha,\gamma}$ (without
offsets) to it to obtain point estimates for what is effectively
$(-\log(\nverts)\alpha^\star+\alpha_0,
\log(\nverts)\gamma^\star+\gamma_0)$. The calculations were
performed using the \texttt{ergm} package (\citepa{HaHu14e}; \citepa{HuHa08e}) for
the R computing environment \citep{R13r}.
%
\begin{figure*}

\includegraphics{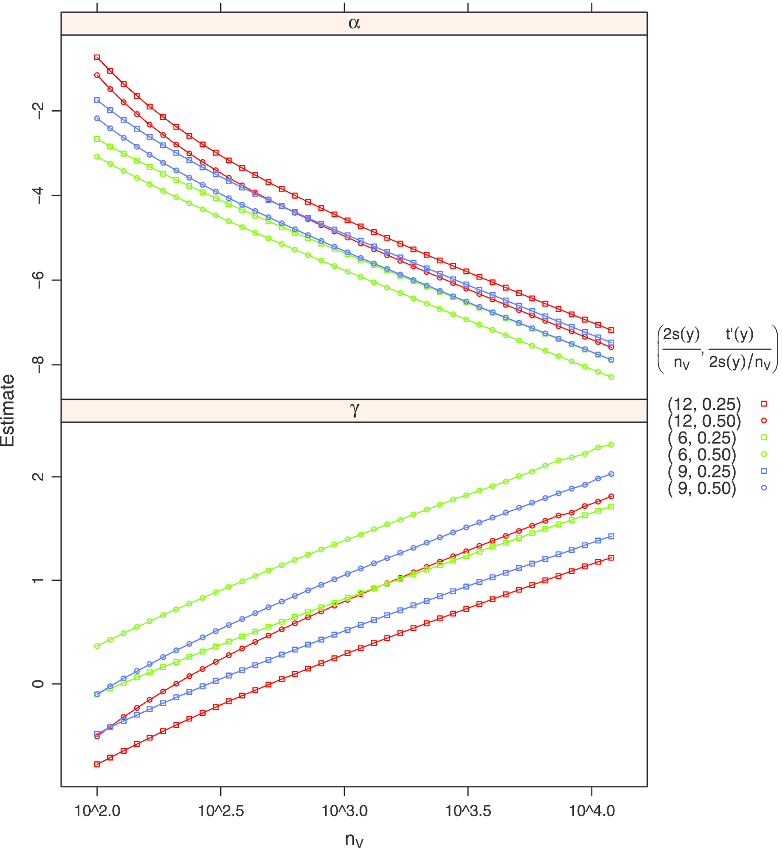}

\caption{Maximum likelihood estimates from fitting
$\pr_{\alpha,\gamma}$ to networks with a variety of sizes, densities
(distinguished by color) and levels of transitivity (distinguished by
plotting symbol). Note that the horizontal axis is plotted on the
logarithmic scale.}\label{fig:triadic-ests}
\end{figure*}

We show the results in Figure~\ref{fig:triadic-ests}. Our
intuition seems to be confirmed, to the extent that our network
sizes are sufficiently large to confirm or disconfirm it. The trend
in both parameter estimates appears to become more linear (in
$\log\nverts$) as $\nverts$ increases, suggesting that unique
$\alpha^\star$ and $\gamma^\star$ exist. For $\hat{\alpha}$, the
asymptotic slope (i.e., $-\alpha^\star$) is very close to $-1$
regardless of the mean degree and the amount of transitivity, and
for $\hat{\alpha}$, the slope (i.e., $\gamma^\star$) decreases as
$\log\nverts$ increases, though it does not quite obtain the exact
value of $1/3$ for the network sizes considered. [Considering
only $\nverts>5\mbox{,}000$, $2\nedges(y)/\nverts=6$, and
$t'(y)/(2\nedges(y)/\nverts)=1/4$---the fastest-converging
configuration---gave the slope of 0.35.]

Notably, even though given a particular value of the
sufficient statistic $(\nedges(y),\ntrties(y))$, the natural
parameters $(\alpha,\gamma)$ would be determined exactly, we have to use
Monte Carlo MLE \citep{HuHa06i} to estimate them, so there is some
noise in the
point estimates.

Overall, it appears that the coefficients of sparser networks
with weaker transitivity tend to approach linearity faster. Thus, we
performed a follow-up simulation study, this one with mean degree 2,
transitivity proportion $1/8$ and 40 values of $\nverts$ between
10\mbox{,}000 and 40\mbox{,}000, logarithmically spaced.

Based on all of the values of $\nverts$ considered,
$\hat{\alpha}^\star= 1.00037$ [95\% CI: $(1.00029,1.00044)$] and
$\hat{\gamma}^\star= 0.3377$ [95\% CI: $(0.3369,0.3386)$], closer to
the theoretical values of $1$ and $1/3$ than the smaller network
sizes. The confidence intervals do not include the theoretical
values, but we would not expect the asymptotic values to be attained
for any finite network size. Indeed, there is evidence of
nonlinearity in that range [$P$-value of predictor $\log(\nverts)^2$
term is
$<$0.0001 for the $\alpha^\star$ response and $0.04$ for the
$\gamma^\star$ response, with negative coefficient for both].
Furthermore, fitting only the 20 data points with $\nverts>20\mbox{,}000$
produces $(\hat{\alpha}^\star,\hat{\gamma}^\star)=(1.000072,
0.3347)$, and fitting only the 10 data points with $\nverts>29\mbox{,}000$,
$(\hat{\alpha}^\star,\hat{\gamma}^\star)=(1.00030,0.3334)$.

This very strongly suggests meaningful and interpretable
asymptotic behavior for triadic closure ERGM terms as well. In
particular, the asymptotic linearity with a known coefficient
suggests a form of consistency for ``intercepts'' $\alpha_0$
and $\gamma_0$, as it is they that control the asymptotic mean
degree and per-vertex amount of triadic closure in
\eqref{eq:bern.w.trt}.

To relate this to the notion of effective sample size $\neff$
used earlier, defined through the
scaling of the information matrix $\cI(\alpha,\gamma)$,
we simulated the sufficient statistics from the
above-described fits. For an exponential family, the
variance--covariance matrix of sufficient statistics under the MLE
approximates the information matrix (\citepa{HuHa06i}, equation 3.5, e.g.). We find that the entries of
$\hat{\cI}(\hat{\alpha},\hat{\gamma})/\nverts=\Var_{\hat
{\alpha},\hat{\gamma}} ( [\nedges(Y),\ntrties(Y)
] )/\nverts$
do not exhibit any trend at all as a function of $\nverts$, for
fixed mean degree and per-vertex transitivity. [In particular, for a
linear trend, $P$-values are 0.31, 0.49 and 0.41 for
$\Var(\nedges(Y))/\nverts$, $\Var(\ntrties(Y))/\nverts$ and
$\Cov(\nedges(Y),\ntrties(Y))/\nverts$, respectively. Exploratory
plots do not show any pattern, except for greater variability in
estimates of variance for higher $\nverts$.] This strongly suggests
that the asymptotics of the model \eqref{eq:bern.w.trt} have
an effective sample size $\neff$ of $\nverts$ as well.

\section{Coverage of Wald confidence intervals}\label{sec:coverage}

Our asymptotic arguments in Section~\ref{sec:results} were developed
primarily for the purpose of establishing the scaling associated with
the asymptotic variance, so as to provide insight into the question of
effective sample size---our main focus here. However, the
asymptotically normal distributions we have derived are of no little
independent interest themselves, as they serve as a foundation for
doing formal inference on the model parameters in practice. By way of
illustration, here we explore their use for constructing confidence
intervals, particularly those based on Theorem~\ref{thm:mutuality}:
under a model $p^\ddag_{\alpha,\beta}$, the Wald confidence
intervals using plug-in estimators for the standard errors are $\hat
{\alpha}^\ddag\pm z^*_{(1-\mathrm{CL})/2} \sqrt{e^{-\hat{\alpha
}^\ddag}/\nverts}$ for $\alpha$ and $\hat{\beta}^\ddag\pm
z^*_{(1-\mathrm{CL})/2} \sqrt{e^{-\hat{\alpha}^\ddag}(4+2e^{-\hat
{\alpha}^\ddag-\hat{\beta}^\ddag})/\nverts}$ for $\beta$.

Because our asymptotics are in $\nverts=|V|$, we examine a variety of
network sizes. The desired asymptotic properties of the network are
expressed in terms of the per-capita mean value parameters---$\E
^\ddag_{\alpha, \beta}[\nedges(Y)/\nverts]$ and $\E^\ddag
_{\alpha,\beta}[2\nmutual(Y)/\nverts]$. We study two configurations:
\renewcommand\thelonglist{(\arabic{longlist})}
\renewcommand\labellonglist{\thelonglist}%
\begin{longlist}[(1)]
\item\label{enum:s1m025} $(\E^\ddag_{\alpha, \beta}[\nedges
(Y)/\nverts], \E^\ddag_{\alpha,\beta}[\nmutual(Y)/\nverts
])=(1,0.25)$\linebreak[4]   and
\item\label{enum:s1m040} $(\E^\ddag_{\alpha, \beta}[\nedges
(Y)/\nverts], \E^\ddag_{\alpha,\beta}[\nmutual(Y)/\nverts])=(1,0.40)$.
\end{longlist}
In other words, the expected mean outdegree is set to $1$, and
expected numbers of \emph{out-ties} that are reciprocated are
$0.25\cdot
2=0.5$ and $0.40\cdot 2=0.8$ per vertex, respectively. These represent two
levels of mutuality, though note that even \ref{enum:s1m025}
represents substantial mutuality, especially for larger networks.

For each $\nverts=10,15,20,\ldots,200$, we estimate the natural
parameters of the model $p^\ddag_{\alpha,\beta}$ corresponding to
the desired mean value parameters,
and then simulate 100\mbox{,}000 networks from each configuration, evaluating
the MLE and constructing a Wald confidence interval at each level of
the customary 80\%, 90\%, 95\% and 99\%, for $\alpha$ and
for $\beta$ (individually), checking the coverage.

For some of the smaller sample sizes, the simulated network
statistics for some realizations were not in the interior of their
convex hull (Barndorff-\break Nielsen \citeyear{Ba78i}, Theorem~9.13, page 151). That is, their
values were
the maximal or minimal possible: $\nedges(y)=0$ (empty graph),
$\nedges(y)=\nverts(\nverts-1)$ (complete graph), $\nmutual(y)=0$
(no ties reciprocated), and/or $\nmutual(y)=\nedges(y)/2$ (every
extant tie
reciprocated). For those, the MLE did not exist. [For
\ref{enum:s1m025}, the fraction was 8.2\% for $\nverts=10$ and none
of the 100\mbox{,}000 realizations had no MLE for $\nverts\ge55$. For
\ref{enum:s1m040}, it was 14.2\% for $\nverts=10$ and none
of the realizations had no MLE for $\nverts\ge65$.]

Our results are conditional on the MLE existing. From the frequentist
perspective, one might argue that if the MLE did not exist for a real
data set, we would not have reported that type of confidence interval,
so it should be excluded from the simulation as well.

\begin{figure*}[b]

\includegraphics{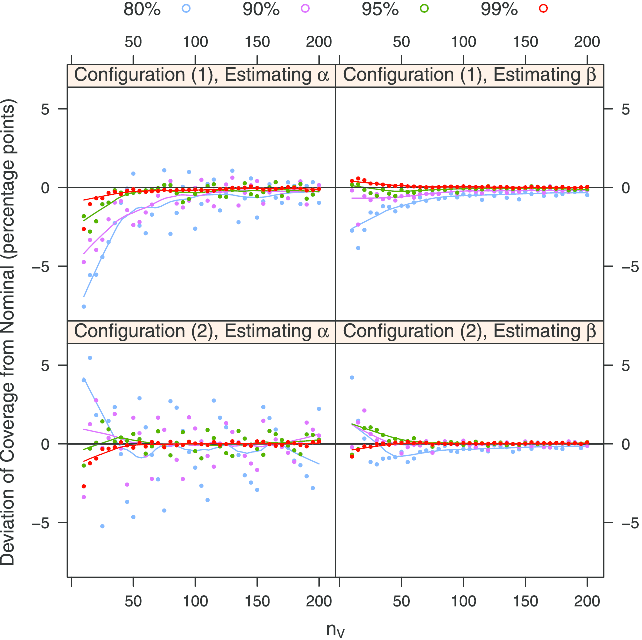}

\caption{Scatterplot of differences between
simulated coverage and nominal coverage for
the two configurations studied, as a function of network size
$\nverts$. Color denotes the nominal coverage levels, and smoothing
lines have been added. Note that the
differences are differences in percentage points ($\mathrm{simulated\ \%
}-\mathrm{nominal\ \%}$), not percent differences
($\frac{\mathrm{simulated\ \%}-\mathrm{nominal\ \%}}{\mathrm{nominal\ \%
}}\cdot 100\%$).}
\label{fig:coverage}
\end{figure*}

We report coverages for selected network sizes in
Table~\ref{tab:coverage} and provide a visualization in
Figure~\ref{fig:coverage}. Overall, the 80\% coverage appears to be
varied---and not very conservative---while higher levels of
confidence appear to be more
consistently conservative, particularly for estimates of $\beta$.
Coverage for $\alpha$ appears to oscillate as a function of network
size. This is particularly noticeable for the lower confidence levels
and stronger mutuality \ref{enum:s1m040}. Tendency of a confidence
interval for a binomial proportion to oscillate around the nominal
level is a known phenomenon (\citepa{brown2001interval}, \citeyear{brown2002confidence}, and
others), though it is
interesting to note that it appears to be more prominent for the
density, rather than mutuality, parameter and that it appears to be
stronger for stronger mutuality.
%
\begin{table*}[t]
\tablewidth=\textwidth
\tabcolsep=0pt
\caption{Simulated Theorem \protect\ref{thm:mutuality}
confidence interval coverage levels for selected network sizes and two
levels of reciprocity: lower \protect\ref{enum:s1m025} and higher \protect\ref
{enum:s1m040}}
\label{tab:coverage}
\begin{tabular*}{\textwidth}{@{\extracolsep{\fill}}lrcccccccc@{}}
\hline
& & \multicolumn{8}{c}{\textbf{Coverage}} \\
\ccline{3-10}
& & \multicolumn{2}{c}{\textbf{80.0\%}} &  \multicolumn{2}{c}{\textbf{90.0\%}} &
\multicolumn{2}{c}{\textbf{95.0\%}} &  \multicolumn{2}{c}{\textbf{99.0\%}}\\
 \ccline{3-4,5-6,7-8,9-10}
& $\bolds{\nverts}$ & $\bolds{\alpha}$ & $\bolds{\beta}$ &  $\bolds{\alpha}$ & $\bolds{\beta}$  &
$\bolds{\alpha}$ & $\bolds{\beta}$ &  $\bolds{\alpha}$ & $\bolds{\beta}$ \\
\hline
\ref{enum:s1m025} & 10 & 72.4\% & 77.3\% &85.3\% & 89.8\% &93.2\%
& 95.2\% &96.4\% & 99.4\% \\
& 20 & 74.5\% & 77.3\% & 86.0\% & 89.4\% &92.9\% & 94.9\% &98.3\% &
99.5\% \\
& 50 & 80.9\% & 78.8\% & 87.6\% & 89.4\% &94.7\% & 94.8\% &98.9\% &
99.2\% \\
& 100 & 77.4\% & 79.6\% & 90.0\% & 90.0\% &94.6\% & 94.9\% &98.9\% &
99.1\% \\
& 200 & 79.0\% & 79.5\% & 90.1\% & 89.8\% &94.9\% & 94.9\% &98.9\% &
99.0\% \\
[3pt]
\ref{enum:s1m040} & 10 & 84.0\% & 84.2\% &86.6\% & 89.8\% &93.6\%
& 94.3\% &96.3\% & 98.2\% \\
& 20 & 81.8\% & 80.3\% & 92.8\% & 92.1\% &95.1\% & 96.0\% &98.1\% &
98.8\% \\
& 50 & 75.3\% & 79.5\% & 91.7\% & 89.4\% &95.6\% & 95.1\% &98.8\% &
99.0\% \\
& 100 & 78.5\% & 79.7\% & 91.0\% & 90.2\% &94.5\% & 94.9\% &99.0\% &
99.1\% \\
& 200 & 82.2\% & 79.9\% & 90.5\% & 89.9\% &95.3\% & 95.1\% &99.2\% &
99.1\% \\
\hline
\end{tabular*}
\end{table*}
%

\section{Example: Food-sharing networks in Lamalera}\label{sec:example}

While the results of Section~\ref{sec:results} are important in
establishing how\vadjust{\goodbreak} closely the question of effective sample size in
network modeling is tied to the structural property of (non)sparseness
expected of the networks modeled, there remains the important practical
question of establishing in applications just which model (i.e., sparse
or nonsparse) is most appropriate. While a full and detailed study of
this question is beyond the scope of this work, we present here an
initial exploration.

\begin{figure*}

\includegraphics{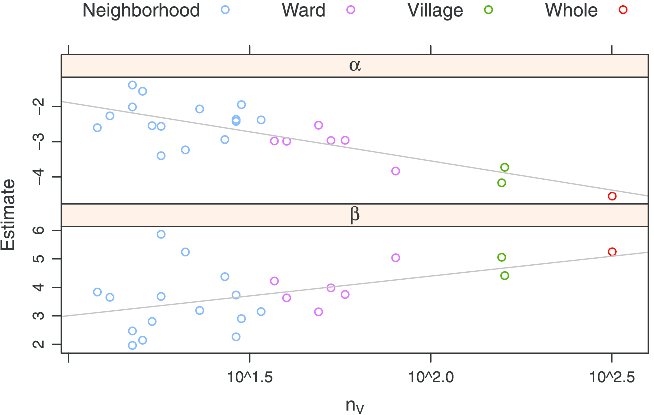}

\caption{Maximum likelihood estimates from fitting
$\pr_{\alpha,\beta}$ to each subdivision of the Lamalera
food-sharing network. Note that the horizontal axis is plotted on the
logarithmic scale. Colors indicate subdivision type. The least-squares
coefficients from regressing $\hat{\alpha}$ and $\hat{\beta}$ on
$\log\nverts$ are $-0.72$ and $+0.60$, respectively.}
\label{fig:lamalera-ests}
\end{figure*}

Note that, in exploring this question, we face a problem similar to that
pointed out by \citet{krivitsky2011adjusting}: it requires a
collection of
closed networks of a variety of sizes yet substantively similar social
structure. Furthermore, our results are limited to modeling density
and reciprocity, so the networks should be well approximated by this
model. Here, we use data collected by \citet{nolin2010food}, in which
each of 317 households in Lamalera, Indonesia was asked to list the
households to whom they have given and households from whom they have
received food in the preceding season. Lamalera is split,
administratively, into two villages, which are further subdivided into
wards, and then into neighborhoods. \citet{nolin2010food} fit several ERGMs
to the network, finding that distance between households had a
significant effect on the propensity to share, as did kinship between
members of the households involved. Nolin also
found a significant positive mutuality effect.

In our study, we make use of the geographic effect by constructing a
series of 24 overlapping subnetworks, consisting of Lamalera itself,
its 2 constituent villages, 6 wards and 15 neighborhoods, with
network sizes ranging from 12 to 317. We then fit the baseline model
$\pr_{\alpha,\beta}$ to each network. If $\pr_{\alpha,\beta}$ is the
most realistic asymptotic regime for these data, we\vspace*{1pt} would expect
estimates $\hat{\alpha}$ and $\hat{\beta}$ to have no relationship to
$\log\nverts$ for the corresponding network. If
$\pr^{\dag}_{\alpha,\beta}$ is the most realistic, we would expect no
relationship between $\log\nverts$ and $\hat{\beta}$, but an
approximately linear relationship with $\hat{\alpha}$, with slope
around $-1$. Last, if $\pr^{\ddag}_{\alpha,\beta}$ is the most
realistic, we would expect the slope of the relationship between $\log
\nverts$ and $\hat{\alpha}$ to be around $-1$ and between $\log
\nverts$ and $\hat{\beta}$ to be around $+1$.

The estimated coefficients and the slopes are given in
Figure~\ref{fig:lamalera-ests}. The results are suggestive. The
relationship between $\hat{\alpha}$ and $\log\nverts$ is clearly
negative, while the relationship between $\hat{\beta}$ and $\log
\nverts$ is clearly positive, and the magnitudes of both slopes are closer
to 1 than to 0 (although both are far from equaling 1). Overlap between
the subnetworks induces dependence among the coefficients, so it is not
possible to formally test or estimate how significant this difference
is. Nevertheless, the preponderance of evidence is that $\pr^{\ddag
}_{\alpha,\beta}$ is the best of the three considered. That is, a
sparse model that does not enforce sparsity on reciprocating ties
appears to be preferable here.

\begin{figure*}

\includegraphics{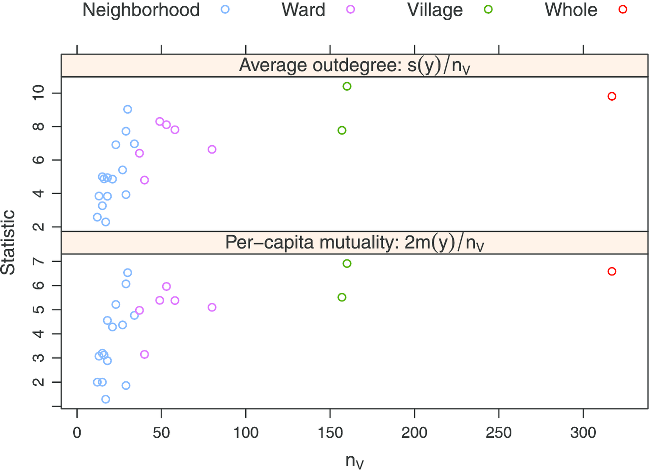}

\caption{Per-capita network statistics as a
function of $\nverts$. Colors indicate subdivision type. Note that
the larger subdivisions have more within-subdivision ties.}
\label{fig:lamalera-stats}
\end{figure*}

A possible explanation for why the magnitudes of the slopes are
substantially less than 1 is that both the argument of
\citet{krivitsky2011adjusting} and our argument in
Section~\ref{sec:bernoulli-reciprocity} rely on the assumption that
the network is closed---no relationships of interest are to or from
vertices outside of the network observed---or, at least, that the
stable mean degree and
per-capita reciprocity are for the ties within it. However, while there
is likely to be very little food
sharing out of or into Lamalera, and relatively little between the two
villages it comprises (7\% of all food-sharing ties in the network are
between villages), there is more sharing between the wards (28\% are
between wards), and even more between neighborhoods (44\%). Thus, the
closed-network assumption is violated. (The respective
between-subdivision percentages for reciprocated ties are 6\%, 22\%
and 39\%.) When each of the subdivisions of the network is considered
in isolation, these ties are lost, so the smaller subdivisions appear,
to the model, to have smaller mean degree and per-capita
mutuality. (See Figure~\ref{fig:lamalera-stats}.) This, in turn, means
that smaller subdivisions have a decreased $\hat{\alpha}$ (increasing
the slope for it in Figure~\ref{fig:lamalera-ests}) and, because mutual
ties suffer less of this ``attrition'' than ties do overall, the
$\hat{\beta}$, after adjusting for the decreased $\hat{\alpha}$, is
increased for smaller networks, thus reducing the slope for
$\hat{\beta}$ in Figure~\ref{fig:lamalera-ests}. It is not unlikely that
this pattern will hold in any network with an unobserved spatial
structure, whose subnetworks of interest are contiguous regions in this space.
%

\section{Discussion}\label{sec:discussion}
Unlike conventional data, network data typically do not have an
unambiguous notion of sample size. The theoretical developments and
the examples we have presented show that the effective sample size
$\neff$
associated with a network depends strongly on the model assumed for
how the network scales. In particular, in the case of reciprocity,
whether or not the model for scaling takes into account the notion of
the preexisting relationship affects whether reciprocity is even
meaningful for large networks. In the case of triadic effects, a
similar notion---along with the intuition that as the network size
changes each individual's view of triadic closure should not---implies a specific scaling regime which, in turn, implies a specific
notion of the effective sample size.

The models we study here are relatively simple examples
of network models.
However, with reciprocity, our work includes an important aspect that
already allows us a glimpse beyond the
treatments of, say, \citet{chatterjee2011random} and \citet
{rinaldo2011maximum},
for so-called beta models, where the dependency induced here by
reciprocity is absent. In addition, the results for reciprocity
suggest that the effective modeling of triadic (e.g., friend of a
friend of a friend) effects in a manner indexed to network size
requires a more complex treatment yet. However, our
simulation shows, perhaps somewhat surprisingly, that if triadic
closure is considered on a per-capita basis, effective sample size
ultimately behaves similarly to the way it does in the simpler cases.

We note that asymptotic theory supporting
methods for the construction of confidence intervals for network parameters
is only beginning to emerge. The most traction appears to have been
gained in
the context of stochastic block models (e.g.,
Bickel and Chen, \citeyear{bickel2009nonparametric};
Choi, Wolfe and Airoldi, \citeyear{choi2012stochastic};
Celisse, Daudin and Pierre, \citeyear{celisse2012consistency} and
Rohe, Chatterjee and Yu, \citeyear{rohe2011spectral}), although
progress is beginning to be had with exponential random graph models as well
(e.g.,
\citepa{chatterjee2011random};
\citepa{chatterjee2013estimating} and
Rinaldo, Petrovi{\'c} and Fienberg, \citeyear{rinaldo2011maximum}). Most of
these works present consistency results for maximum likelihood and
related estimators, with the exception of \citet
{bickel2009nonparametric}, which also includes results on asymptotic
normality of estimators. Our work contributes to this important but
nascent area with both our theoretical developments and our simulation
studies. In particular, the asymptotic regime of $\pr_{\alpha,\gamma
}$ is one that neither appears to become degenerate nor approaches
Erd\H{o}s--R\'enyi.

The lack of an established understanding of the distributional
properties of parameter estimates in commonly used network models is
particularly unfortunate given that a number of software packages now
allow for the easy computation of such estimates. For example, packages
for computing estimates of parameters in fairly general formulations of
exponential random graph models routinely
report both estimates and, ostensibly, standard errors, where the
latter are based on standard arguments for exponential families.
Unfortunately, practitioners do not always seem to be aware that the
use of these standard errors for constructing normal-theory confidence
intervals and tests is lacking fully formal justification. From
that perspective, our work appears to be one of the first to begin
laying the necessary
theoretical foundation to justify practical confidence interval
procedures in exponential
random graph models. See \citet{haberman1981exponential} for another
contribution in this direction, proposed
as part of the discussion of the original paper of \citet
{holland1981exponential}.

In order to successfully build upon our work, and extend our results
to more sophisticated instances of exponential random graph models,
certain technical
challenges must be overcome. First, we note that our notion of
effective sample size
is tied directly to the asymptotic behavior of the Fisher information matrix
of our model [denoted $\mathcal{I}(\theta)$ in the supplemental
article \citep{Suppl}].
Given that exponential random graph models are, by definition,
of exponential family form, this information matrix is in principle given
by the matrix of partial second derivatives of the cumulant generating function
[denoted $\psi$ in the supplemental article, so that
$\mathcal{I}=\partial^2 \psi(\theta)/\partial\theta\,\partial
\theta^\top$].
Due to the use of independent dyads in our theoretical work
(i.e., our models are variations on Bernoulli models), the corresponding
likelihoods factor over dyads, and hence the information matrices
are simply proportional to powers of $\nverts$ (i.e., linear or quadratic).
This is in analogy to the canonical setting of independent and identically
distributed observations. However, in more general settings beyond the
case of independent dyads---including even the models with triadic
effects we studied in simulation---the likelihood cannot be expected to
factor in such a simple manner. Hence, the analysis of the Fisher information
promises to be decidedly more subtle. In fact, there appears to be almost
no work to date studying this matrix in any detail. To the best of our
knowledge,
the only such work is the recent manuscript by \citet{pu2013learning},
introducing a deterministic approach to approximating this matrix (stochastic
approximations may, of course, be produced using MCMC) based on a lower bound
of the cumulant generating function. This bound, however, has only an implicit
representation.

Second, in the case of more general exponential random graph models
than those studied here, there will be a need for a correspondingly more
sophisticated central limit theorem, in order to produce results on
asymptotic normality analogous to those we present for the simpler
models we
study. Even for our models, the tool we used was somewhat
nonstandard, in that we required a double-array central limit theorem.
The more general case will require a central limit theorem capable of
handling the nontrivial global dependencies induced by effects even as
seemingly simple as triadic closure or the like. Progress on the first point
above is a likely prerequisite to understanding the nature of these dependencies
sufficiently well to know just what sort of central limit theorem is required.

Finally, there is, as always with exponential random graph models,
the issue of instability and degeneracy that must be kept in
mind (e.g., \citepa{handcock.2003} and \citepa{chatterjee2011random}).
It has been discovered only relatively recently that substantial
care must be taken in specifying network effects in exponential
random graph models. Without such care, it is possible to produce
models for which the corresponding distributions turn out to be
near-degenerate and, in turn, the estimation of parameters highly
unstable. \citet{schweinberger2011instability} has recently shed
important light on this issue, showing that instability and degeneracy
are related to the scaling of the linear term in exponential family
distributions generally and, more specifically, in exponential random
graph models. These scaling results can be expected to have implications
on the role that scaling necessarily plays in the types of calculations
we have presented here.


\section*{Acknowledgments}
Pavel N. Krivitsky was supported in part by ONR
Award N000140811015, NIH Award
1R01HD068395-01, and Portuguese Foundation for Science
and Technology Ci\^encia 2009 Program and also wishes
to thank David Nolin for the Lamalera food-sharing
network data set. Eric D. Kolaczyk was supported
in part by ONR Award N000140910654. This work
was begun during the 2010--2011
Program on Complex Networks at SAMSI.

\begin{supplement}[id=suppA]
\stitle{Supplement to ``On the Question of Effective Sample Size in Network Modeling:
An Asymptotic Inquiry''}
\slink[doi]{10.1214/14-STS502SUPP} 
\sdatatype{.pdf}
\sfilename{sts502\_supp.pdf}
\sdescription{This document contains proofs of the results reported in the body of the article.}
\end{supplement}


\end{document}